\documentstyle[12pt]{article}
\input{prepictex}
\input{pictex}
\input{postpictex}

\begin{document}
\vsize=24.0truecm
\hfill
\vskip 3.8truecm
\thispagestyle{empty}
\def\thefootnote{}
\newcommand{\R}{\mbox{$I\!\! R$}}
\centerline{\large A discrete form of the Beckman-Quarles theorem for}
\centerline{\large two-dimensional strictly convex normed spaces}
\vskip 0.8truecm
\centerline{\large Apoloniusz Tyszka}
\vskip 0.8truecm
\centerline{version of October 8, 2000}
\footnotetext{
\footnotesize
\par
\noindent
submitted to {\it Aequationes Mathematicae}
}
\hfill
\par
\begin{abstract}
Let $X$ be a real normed vector space and ${\rm dim} X \ge 2$.
Let $\rho > 0$ be a fixed real number. We prove that
if $x,y\in X$ and $||x-y||/\rho$ is a rational number
then there exists a finite set $\{x,y\} \subseteq S_{xy}\subseteq X$
with the following property: for each strictly convex $Y$ of
dimension $2$ each map from $S_{xy}$ to $Y$ preserving the distance
$\rho$ preserves the distance between $x$ and $y$. It implies that
each map from $X$ to $Y$ that preserves the distance $\rho$ is an
isometry.
\end{abstract}
\vskip 0.1truecm
\par
Let {\bf Q} denote the field of rational numbers.
All vector spaces mentioned in this article are assumed to be real.
A normed vector space $E$ is called {\sl strictly convex} ([5]),
if for each pair $a$, $b$ of nonzero elements in $E$ such that
$||a+b||=||a||+||b||$, it follows that $a=\gamma b$ for some $\gamma>0$.
It is known~([15]) that two-dimensional strictly convex normed spaces
satisfy the following condition~$(\ast)$:
\begin{description}
\item{$(\ast)$}
for any $a \neq b$ on line $L$ and any $c,d$ on the same
side of $L$, if $||a-c||=||a-d||$ and $||b-c||=||b-d||$, then $c=d$.
\end{description}
Conversely ([15]), for any two-dimensional normed space
the condition $(\ast)$ implies that the space is strictly convex.
\par
The classical Beckman-Quarles theorem states that any map from
${\R}^{n}$ to ${\R}^{n}$ ($2\leq n<\infty $) preserving unit
distance is an isometry, see [1], [2] and~[6].
Various unanswered questions and counterexamples concerning
the Beckman-Quarles theorem and isometries are discussed
by Ciesielski and Rassias [4]. For more open problems and
new results on isometric mappings the reader is referred to [7]-[13].
The Theorem below may be viewed as a discrete form of the
Beckman-Quarles theorem for two-dimensional strictly convex
normed spaces.
\vskip 0.3truecm
\par
{\bf Theorem}.
Let $X$ and $Y$ be normed vector spaces such that
${\rm dim} X \ge {\rm dim} Y = 2$ and $Y$ is strictly convex.
Let $\rho>0$ be a fixed real number.
\vskip 0.1truecm
\par
\noindent
{\bf 1.} If $x,y\in X$ and $||x-y||/\rho$ is a rational number
then there exists a finite set $S_{xy} \subseteq X$ containing
$x$ and $y$ such that each injective map $f:S_{xy} \rightarrow Y$
preserving the distance $\rho$ preserves the distance between
$x$ and $y$.
\vskip 0.1truecm
\par
\noindent
{\bf 2.} If $x,y\in X$ and $\varepsilon >0$ then there exists a
finite set $T_{xy}(\varepsilon) \subseteq X$ containing $x$ and
$y$ such that each injective map
$f:T_{xy}(\varepsilon )\rightarrow Y$ preserving the
distance~$\rho$ preserves the distance between $x$ and $y$ to
within $\varepsilon $ i.e.
\begin{displaymath}
|||f(x)-f(y)||-||x-y||| \leq \varepsilon.
\end{displaymath}
\vskip 0.1truecm
\par
\noindent
{\bf 3.} Let $X={\R}^n$ ($2 \le n < \infty$) be equipped with euclidean
norm. Then the assumption of injectivity is unnecessary in items
1 and 2.
\vskip 0.1truecm
\par
\noindent
{\bf 4.} More generally (cf. item 3), for each normed space $X$
the assumption of injectivity is unnecessary in items 1 and 2.
\vskip 0.3truecm
\par
\noindent
{\bf Proof of item 1.}
Let $D$ denote the set of all non-negative numbers $d$ with the
following property $(\ast\ast)$:
\begin{description}
\item{$(\ast\ast)$}
if $x,y\in X$ and $||x-y||=d$ then there exists a finite set
$S_{xy}\subseteq X$ such that $x,y\in S_{xy}$ and any injective
map $f:S_{xy}\rightarrow Y$ that preserves the distance $\rho$
also preserves the distance between $x$ and $y$.
\end{description}
\vskip 0.1truecm
\par
\noindent
Obviously $0,\rho \in D$. We first prove that if $d \in D$, then
$2 \cdot d \in D$. Assume that $d \in D$, $d>0$, $x,y \in X$,
$||x-y||=2 \cdot d$. Using the notation of Figure 1
\vskip 0.3truecm

\centerline{
\beginpicture
\normalsize
\setcoordinatesystem units <1mm, 1mm>
\setplotarea x from -3 to 77, y from -6 to 23
\setplotsymbol({\xxpt\rm.})
\plot 0 0 37 0 74 0 62 20 25 20 37 0 62 20 25 20 0 0 /
\put {$x$} at 0 -3
\put {$z$} at 37 -3
\put {$y$} at 74 -3
\put {$x_1$} at 62 23
\put {$y_1$} at 25 23
\put {$d$} at 18.5 2
\put {$d$} at 55.5 2
\put {$d$} at 34 10
\put {$d$} at 8 10
\put {$d$} at 43.5 22
\put {$d$} at 71 10
\put {$d$} at 45 10
\endpicture}
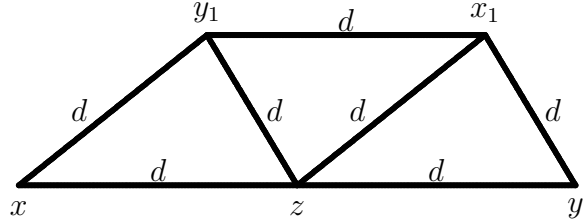
\centerline{Figure 1}
\centerline{$||x-y||=2 \cdot d$}
\centerline{$z:=\frac{x+y}{2}$}
\centerline{$||x-z||=||x-y_1||=||z-y_1||=d$}
\centerline{$x_1:=y_1+(z-x)$}
\par
\noindent
we show that
\\
\centerline{$S_{xy}:=S_{xz} \cup S_{zy} \cup S_{y_1x_1} \cup 
S_{xy_1} \cup S_{zx_1} \cup S_{zy_1} \cup S_{yx_1}$}
\\
satisfies the condition $(\ast\ast)$.
Let an injective $f:S_{xy} \rightarrow Y$ preserves the distance
$\rho$. By the injectivity of $f$:
$f(x) \neq f(x_1)$ and $f(y) \neq f(y_1)$. According to $(\ast)$:
$f(y_1)-f(x_1)=f(x)-f(z)$ and $f(y_1)-f(x_1)=f(z)-f(y)$.
Hence $f(x)-f(z)=f(z)-f(y)$.
Therefore $||f(x)-f(y)||=||2(f(x)-f(z))||=2 \cdot ||f(x)-f(z)||=
2 \cdot ||x-z||=2 \cdot d =||x-y||$.
\vskip 0.2truecm
From Figure 2, the previous step and the property that defines strictly
convex normed spaces it is clear that if $d \in D$, then all distances
$k \cdot d$ ($k$ a positive integer) belong to $D$.

\centerline{
\beginpicture
\normalsize
\setcoordinatesystem units <1mm, 1mm>
\setplotarea x from -3 to 67, y from -3 to 3
\setplotsymbol({\xxpt\rm.})
\plot 0 0 20 0 40 0 64 0 /
\plot 76 0 80 0 100 0 /
\put {$\bullet$} at 0 0
\put {$\bullet$} at 20 0
\put {$\bullet$} at 40 0
\put {$\bullet$} at 60 0
\put {$\bullet$} at 80 0
\put {$\bullet$} at 100 0
\put {$d$} at 10 2
\put {$d$} at 30 2
\put {$d$} at 50 2
\put {$d$} at 90 2
\begin{scriptsize}
\put {$\bullet$} at 67 0
\put {$\bullet$} at 70 0
\put {$\bullet$} at 73 0
\end{scriptsize}
\put {$x=w_0$} at 0 -3
\put {$w_{1}$} at 20 -3
\put {$w_{2}$} at 40 -3
\put {$w_{3}$} at 60 -3
\put {$w_{k-1}$} at 80 -3
\put {$w_{k}=y$} at 100 -3
\endpicture}
\vskip 0.1truecm
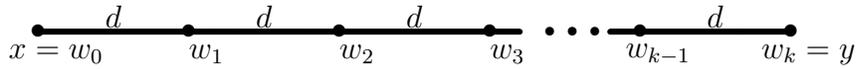
\centerline{Figure 2}
\centerline{$||x-y||=k \cdot d$}
\centerline{$S_{xy}=\bigcup \{S_{ab}:a,b\in \{w_{0},w_{1},...,w_{k}\},
||a-b||=d \vee ||a-b||=2 \cdot d \}$}
\vskip 0.2truecm
\par
From Figure 3, the previous step and the property that defines
strictly convex normed spaces it is clear that if $d \in D$,
then all distances $d/k$ ($k$ a positive integer) belong to $D$.
Hence $D/\rho:=\{d/\rho: d \in D\} \supseteq {\bf Q} \cap [0,\infty)$.
This completes the proof of item 1.
\vskip 0.3truecm

\centerline{
\beginpicture
\normalsize
\setcoordinatesystem units <1mm, 1mm>
\setplotarea x from -3 to 28, y from -18 to 30
\setplotsymbol({\xxpt\rm.})
\plot 10 0 -40 25 -40 -15 10 0 0 5 0 -3 /
\put {$x$} at 0 -5
\put {$y$} at 0 8
\put {$z$} at 12 0
\put {$\tilde{x}$} at -40 -18
\put {$\tilde{y}$} at -40 28
\put {$d$} at -43 5
\put {$d$} at 5 5
\put {$d$} at 5 -4
\put {$(k-1) \cdot d$} at -11 17
\put {$(k-1) \cdot d$} at -11 -11
\endpicture}
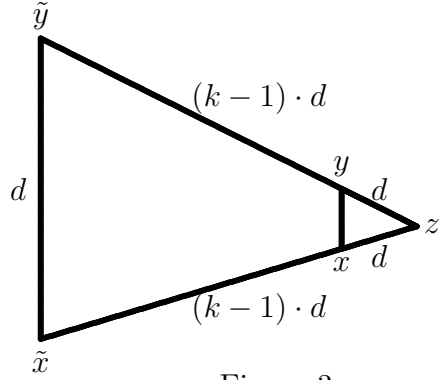
\centerline{Figure 3}
\centerline{$||x-y||=d/k$}
\centerline{$\tilde{x}:=x+(k-1)(x-z)$}
\centerline{$\tilde{y}:=y+(k-1)(y-z)$}
\centerline{$\tilde{x}-\tilde{y}=x-y+(k-1)((x-z)-(y-z))=k(x-y)$}
\centerline{$S_{xy}=S_{\tilde{x}\tilde{y}}
\cup
S_{\tilde{x}x}
\cup
S_{xz}
\cup
S_{\tilde{x}z}
\cup
S_{\tilde{y}y}
\cup
S_{yz}
\cup
S_{\tilde{y}z}$}
\vskip 0.3truecm
\par
\noindent
{\bf Proof of item 2.} It follows from Figure 4.

\centerline{
\beginpicture
\normalsize
\setcoordinatesystem units <1mm, 1mm>
\setplotarea x from -3 to 66, y from -3 to 13
\setplotsymbol({\xxpt\rm.})
\plot 0 0 60 10 63 0 /
\plot 0 0 3 0 /
\plot 6 0 9 0 /
\plot 12 0 15 0 /
\plot 18 0 21 0 /
\plot 24 0 27 0 /
\plot 30 0 33 0 /
\plot 36 0 39 0 /
\plot 42 0 45 0 /
\plot 48 0 51 0 /
\plot 54 0 57 0 /
\plot 60 0 63 0 /
\put {$x$} at 0 -3
\put {$y$} at 63 -3
\put {$z$} at 60 13
\endpicture}
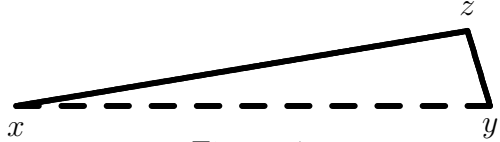
\centerline{Figure 4}
\centerline{$|x-z|/\rho,|z-y|/\rho \in
{\bf Q} \cap [0,\infty)$, $|z-y|\leq \varepsilon /2$}
\centerline{$T_{xy}(\varepsilon )=S_{xz}\cup S_{zy}$}
\par
\noindent
{\bf Proof of item 3.}
In proofs of items 1 and 2 the assumption of injectivity is necessary
only in the first step for distances $2 \cdot d$, $d \in D$.
Let $X=\R^n$ ($2 \le n < \infty$) be equipped with euclidean norm and
$D$ is defined without the assumption of injectivity. Let $d \in D$,
$d>0$. We need to prove that $2 \cdot d \in D$. Let us see at
configuration from Figure 5 below, all segments have the length~$d$.

\centerline{
\beginpicture
\normalsize
\setcoordinatesystem units <1cm, 1cm>
\setplotarea x from -1.00 to 13.00, y from -5 to 6.50
\setplotsymbol({\xxpt\rm.})
\plot 6.00 0.00 0.00 0.00 3.00 5.20 6.00 0.00 9.00 5.20 3.00 5.20 /
\plot 9.38 -0.80 4.38 -4.10 4.00 1.88 9.38 -0.80 9.00 5.20 4.00 1.88 /
\plot 0.00 0.00 4.38 -4.10 /
\plot 6.00 0.00 12.00 0.00 9.00 5.20 6.00 0.00 3.00 5.20 9.00 5.20 /
\plot 2.62 -0.80 7.62 -4.10 8.00 1.88 2.62 -0.80 3.00 5.20 8.00 1.88 /
\plot 12.00 0.00 7.62 -4.10 /
\put {$x$} at -0.30 0.00
\put {$z$} at 6.00 -0.30
\put {$y$} at 12.30 0.00
\put {$x_1$} at 9.00 5.50
\put {$y_1$} at 3.00 5.50
\put {$\tilde{x}$} at 4.38 -4.40
\put {$\tilde{y}$} at 7.62 -4.40
\put {$\tilde{x_1}$} at 8.30 1.88
\put {$\tilde{y_1}$} at 3.70 1.88
\put {$z_x$} at 2.32 -0.80
\put {$z_y$} at 9.68 -0.80
\endpicture}
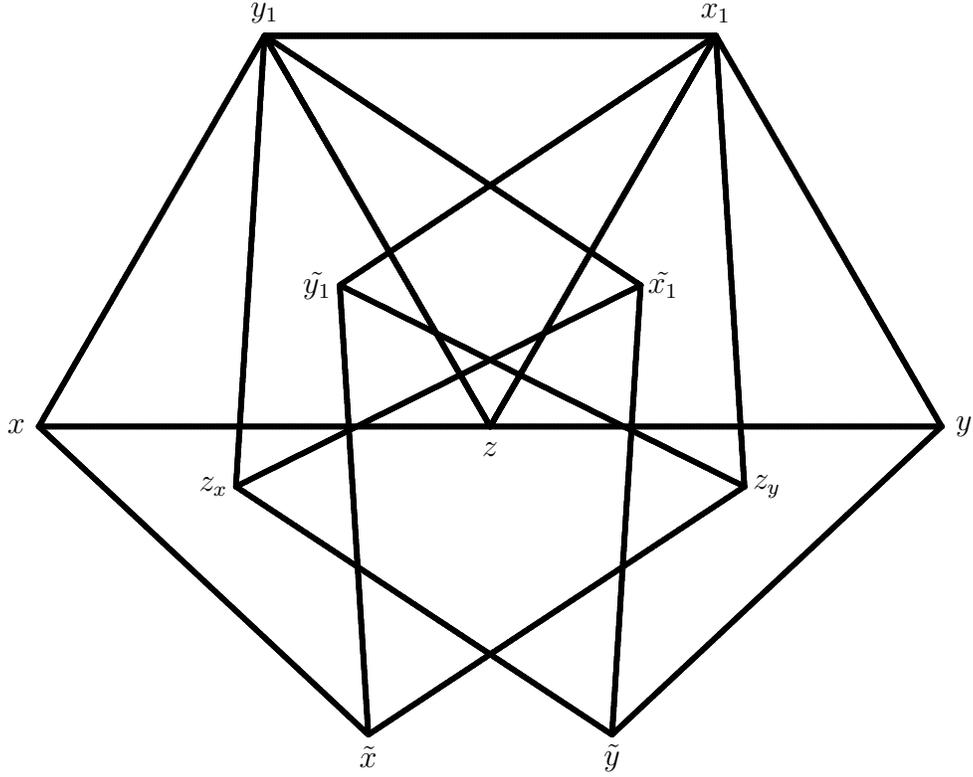
\centerline{Figure 5}
\centerline{$||x-y||=2 \cdot d$}
\centerline{$z:=\frac{x+y}{2}$}
\centerline{$S_{xy}=\bigcup \{S_{ab}:a,b \in
\{x,\tilde{x},x_1,\tilde{x_1},y,\tilde{y},y_1,\tilde{y_1},z,z_x,z_y\},
||a-b||=d\}$}
\vskip 0.3truecm
\par
Assume that $f:S_{xy} \rightarrow Y$ preserves the distance $\rho$.
It is sufficient to prove that $f(x)\neq f(x_1)$ and similarly
$f(y) \neq f(y_1)$. Suppose, on the contrary, that $f(x)=f(x_1)$,
the proof of $f(y) \neq f(y_1)$ is similar. Hence four points:
$f(\tilde x)$, $f(z_y)$, $f(\tilde{y_1})$, $f(x_1)$ have the
distance $d$ from each other. We prove that it is impossible in
two-dimensional strictly convex normed spaces. Suppose, on the
contrary, that $a_1,a_2,a_3,a_4 \in Y$ and $||a_1-a_2||=||a_1-a_3||=
||a_1-a_4||=||a_2-a_3||=||a_2-a_4||=||a_3-a_4||=d>0$.
Let us consider the segment $a_2a_3$. According to $(\ast)$
$a_1$ and $a_4$ lie on the opposite sides of the line $L(a_2,a_3)$
and $a_2-a_1=a_4-a_3$. Let us consider the segment $a_1a_3$.
According to $(\ast)$ $a_2$ and $a_4$ lie on the
opposite sides of the line $L(a_1,a_3)$ and $a_1-a_2=a_4-a_3$.
Hence $a_4-a_3=0$, a contradiction.
This completes the proof of item 3.
\vskip 0.3truecm
\par
\noindent
{\bf Proof of item 4.}
Analogously as in the proof of item 3 it suffices to prove that for
each $x,y \in X$, $x \neq y$ there exist points forming the
configuration from Figure 5 where all segments have the length
$||x-y||/2$. Let us consider $x,y \in X$, $x \neq y$. We choose
two-dimensional subspace $\tilde{X} \subseteq X$ containing $x$ and
$y$.
\par
{\bf First case:} the norm induced on $\tilde{X}$ is strictly convex.
Obviously $\tilde{X}$ is isomorphic to ${\R}^2$ as a linear
space. Let us consider ${\R}^2$ with a strictly convex norm
$||$ $||$. It suffices to prove that for each $a,b \in {\R}^2$
satisfying $||a||=||b||=||a-b||=d>0$ there exist
$\tilde{a},\tilde{b} \in {\R}^2$
satisfying $||\tilde{a}||=||\tilde{b}||=||\tilde{a}-\tilde{b}||=
||(\tilde{a}+\tilde{b})-(a+b)||=d$. We fix $a=(a_x,a_y)$ and
$b=(b_x,b_y)$. Let $S:=\{x \in {\R}^2: ||x||=d\}$. According to
$(\ast)$ for each $u=(u_x,u_y) \in S$ there exists a unique
$h(u)=(h(u)_x,h(u)_y) \in S$ such that $||u-h(u)||=d$ and
\begin{center}
\begin{math}
$$
\det
\left[\begin{array}{cccc}
u_x & u_y \\
h(u)_x & h(u)_y \\
\end{array} \right]
\cdot
\det
\left[\begin{array}{cccc}
a_x & a_y \\
b_x & b_y \\
\end{array} \right]
> 0.
$$
\end{math}
\end{center}
\vskip 0.1truecm
\noindent
Obviously $h(a)=b$.
The mapping $h:S \rightarrow S$ is continuous. For each $u \in S$
$h(-u)=-h(u)$ and
$||u+h(u)||=||2u-(u-h(u))|| \ge |||2u||-||u-h(u)|||=d.$
The following function
\vskip 0.1truecm
\par
\noindent
\centerline 
{$S \ni x \stackrel{g}{\longrightarrow}
||x+h(x)-a-h(a)|| \in [0,\infty)$}
\\
is continuous. We have:
\\
\centerline{$g(a)=0$,}
\centerline{$g(-a)=||-a+h(-a)-a-h(a)||=2\cdot||a+h(a)|| \ge 2 \cdot d$.}
\\
Since $g$ is continuous there exists
$\tilde{a} \in S$ such that $g(\tilde{a})=d$. From this
$\tilde{a}$ and $\tilde{b}:=h(\tilde{a})$ satisfy
$||\tilde{a}||=||\tilde{b}||=||\tilde{a}-\tilde{b}||=
||(\tilde{a}+\tilde{b})-(a+b)||=d$. This completes the proof of
item 4 in the case where the norm induced on $\tilde{X}$ is strictly
convex.
\par
{\bf Second case:} we assume only that $||$ $||$ is a norm on
$\tilde{X}$. The graph $\Gamma$ from Figure 5 (11 vertices, 19 edges)
has the following matrix representation:
\par
\noindent
\begin{center}
\begin{tabular}{|c|c|c|c|c|c|c|c|c|c|c|c|} \hline
$~$ & $v_0$ & $v_1$ & $v_2$ & $v_3$ & $v_4$ & $v_5$ & $v_6$ & $v_7$ & $v_8$ & $v_9$ & $v_{10}$ \\ \hline
$v_0:=x$             & $0$ & $0$ & $1$ & $1$ & $0$ & $0$ & $0$ & $1$ & $0$ & $0$ & $0$ \\ \hline
$v_1:=y$             & $0$ & $0$ & $1$ & $0$ & $1$ & $0$ & $1$ & $0$ & $0$ & $0$ & $0$ \\ \hline
$v_2:=\frac{x+y}{2}$ & $1$ & $1$ & $0$ & $0$ & $1$ & $0$ & $0$ & $1$ & $0$ & $0$ & $0$ \\ \hline
$v_3:=\tilde{x}$     & $1$ & $0$ & $0$ & $0$ & $0$ & $0$ & $0$ & $0$ & $1$ & $0$ & $1$ \\ \hline
$v_4:=x_1$           & $0$ & $1$ & $1$ & $0$ & $0$ & $0$ & $0$ & $1$ & $1$ & $0$ & $1$ \\ \hline
$v_5:=\tilde{x_1}$   & $0$ & $0$ & $0$ & $0$ & $0$ & $0$ & $1$ & $1$ & $0$ & $1$ & $0$ \\ \hline
$v_6:=\tilde{y}$     & $0$ & $1$ & $0$ & $0$ & $0$ & $1$ & $0$ & $0$ & $0$ & $1$ & $0$ \\ \hline
$v_7:=y_1$           & $1$ & $0$ & $1$ & $0$ & $1$ & $1$ & $0$ & $0$ & $0$ & $1$ & $0$ \\ \hline
$v_8:=\tilde{y_1}$   & $0$ & $0$ & $0$ & $1$ & $1$ & $0$ & $0$ & $0$ & $0$ & $0$ & $1$ \\ \hline
$v_9:=z_x$           & $0$ & $0$ & $0$ & $0$ & $0$ & $1$ & $1$ & $1$ & $0$ & $0$ & $0$ \\ \hline
$v_{10}:=z_y$        & $0$ & $0$ & $0$ & $1$ & $1$ & $0$ & $0$ & $0$ & $1$ & $0$ & $0$ \\ \hline
\end{tabular}
\end{center}
\par
\noindent
Let $u_0:=v_0=x$, $u_1:=v_1=y$, $u_2:=v_2=\frac{x+y}{2}$.
We define the following function $\psi$:
\par
\noindent
\centerline{$
\tilde{X}^{8} \ni (u_{3},...,u_{10})
\stackrel{\psi}{\longrightarrow}
(||u_i-u_j||: 0 \le i<j \le 10, (v_i,v_j) \in \Gamma) \in \R^{19}$.}
\par
\noindent
The image of $\psi$ is a closed subset of $\R^{19}$.
For each $\varepsilon>0$ and each bounded $B \subseteq \tilde{X}$
the norm $||$ $||$ may be approximate on $B$ with
$\varepsilon$-accuracy by a strictly convex norm on $\tilde{X}$.
Therefore according to the first case for each
$x,y \in X$, $x \neq y$ and each $\varepsilon>0$ there exist
points forming the configuration from Figure 5 where all segments
have $||$ $||$-lengths belonging to the interval
$(\frac{||x-y||}{2} - \varepsilon,\frac{||x-y||}{2} + \varepsilon)$.
Therefore:
\par
\noindent
\centerline{$(||x-y||/2,...,||x-y||/2) \in \overline{\psi(\tilde{X}^{8})}$
 (the closure of $\psi(\tilde{X}^{8})$).}
\newpage
\par
\noindent
Since $\psi(\tilde{X}^{8})$ is closed we conclude that
\par
\noindent
\centerline{$(||x-y||/2,...,||x-y||/2) \in \psi(\tilde{X}^{8})$.}
This completes the proof of item 4.
\vskip 0.1truecm
\par
{\bf Corollary.}
Let $X$ and $Y$ be normed vector spaces such that
${\rm dim} X \ge {\rm dim} Y = 2$ and $Y$ is strictly convex.
From item 2 of the Theorem follows that an injective map
$f:X \rightarrow Y$ that preserves a fixed distance $\rho>0$ is an
isometry. According to item 4 of the Theorem the assumption
of injectivity is unnecessary in the above statement.
\vskip 0.1truecm
\par
{\bf Remark 1.} The set $S_{xy}$ constructed in the proof does not
depend on $Y$.
\vskip 0.1truecm
{\bf Remark 2.} Instead of injectivity in the Theorem and Corollary
we may assume that
\\
\centerline{$\forall u,v \in {\rm dom}(f) (||u-v||/\rho \in {\bf Q}
\cap (0,\infty) \Rightarrow ||f(u)-f(v)|| \neq ||u-v||/2)$}
\\
It follows from Figure 1.
\vskip 0.1truecm
\par
{\bf Remark 3.}
W. Benz and H. Berens proved ([3], see also [2] and [10]) the
following theorem:
Let $X$ and $Y$ be normed vector spaces such that $Y$ is
strictly convex and such that the dimension of $X$ is at least
$2$. Let $\rho>0$ be a fixed real number and let $N>1$ be a fixed
integer. Suppose that $f :X\rightarrow Y$ is a mapping satisfying:
\\
\centerline{$||a-b||=\rho \Rightarrow ||f(a)-f(b)|| \le \rho$}
\centerline{$||a-b||=N \rho \Rightarrow ||f(a)-f(b)|| \ge N \rho $}
for all $a,b \in X$. Then $f$ is an affine isometry.
\vskip 0.1truecm
\par
{\bf Remark 4}.
A. Tyszka proved ([14]) the following theorem:
if $x,y\in {\R}^{n}$ ($2 \le n < \infty $) and $|x-y|$ is an algebraic
number then there exists a finite set $S_{xy} \subseteq {\R}^{n}$
containing $x$ and $y$ such that each map from $S_{xy}$ to ${\R}^{n}$
preserving unit distance preserves the distance between $x$ and $y$.
\newpage
\centerline{{\bf References}}
\begin{enumerate}
\item F. S. Beckman and D. A. Quarles Jr.,
On isometries of euclidean spaces,
{\it Proc. Amer. Math. Soc.} 4 (1953), 810--815.
\item W. Benz,
{\it Real geometries}, BI Wissenschaftsverlag, Mannheim, 1994.
\item W. Benz and H. Berens, A contribution to a theorem of
Ulam and Mazur, {\it Aequationes Math.} 34 (1987), 61--63.
\item K. Ciesielski and Th. M. Rassias, On some properties of
isometric mappings, {\it Facta Univ. Ser. Math. Inform.}
7 (1992), 107--115.
\item J. A. Clarkson, Uniformly convex spaces,
{\it Trans. Amer. Math. Soc.} 40 (1936), 396--414.
\item U. Everling, Solution of the isometry problem stated
by K. Ciesielski,
{\it Math. Intelligencer}, 10 (1988), No.4, p.47.
\item B. Mielnik and Th. M. Rassias, On the Aleksandrov problem of
conservative distances, {\it Proc. Amer. Math. Soc.}
116 (1992), 1115--1118.
\item Th. M. Rassias, Is a distance one preserving mapping
between metric spaces always an isometry?,
{\it Amer. Math. Monthly} 90 (1983), p.200.
\item Th. M. Rassias, Some remarks on isometric mappings,
{\it Facta Univ. Ser. Math. Inform.} 2 (1987), 49--52.
\item Th. M. Rassias, Properties of isometries and approximate
isometries, {\it in} ''Recent progress in inequalities''
(ed. G. V. Milovanovi\'c), 341--379, Math. Appl. 430,
Kluwer Acad. Publ., Dordrecht, 1998.
\item Th. M. Rassias, Properties of isometric mappings,
{\it J. Math. Anal. Appl.} 235 (1999), 108--121.
\item Th. M. Rassias, Isometries and approximate isometries,
{\it Int. J. Math. Math. Sci.}, to appear.
\item Th. M. Rassias and P. \v{S}emrl, On the Mazur-Ulam theorem and
the Aleksandrov problem for unit distance preserving mappings,
{\it Proc. Amer. Math. Soc.} 118 (1993), 919--925.
\item A. Tyszka, Discrete versions of the Beckman-Quarles theorem,
{\it Aequationes Math.} 59 (2000), 124--133.
\item J. E. Valentine, Some implications of Euclid's proposition 7,
{\it Math. Japon.} 28 (1983), 421--425.
\end{enumerate}
\vskip 0.1truecm
Technical Faculty
\\
Hugo Ko{\l}{\l}\c{a}taj University
\\
Balicka 104, 30-149 Krak\'ow, Poland
\\
E-mail: {\it rttyszka@cyf-kr.edu.pl}
\\
Home page: {\it http://www.cyf-kr.edu.pl/\symbol{126}rttyszka}
\end{document}